\documentclass[letterpaper, 10 pt, conference]{ieeeconf}%
\usepackage{graphics}
\usepackage{epsfig}
\usepackage{cite}
\usepackage{breakurl}
\usepackage{amsmath}
\usepackage{amssymb}
\usepackage{amsfonts}
\usepackage{graphicx}%
\providecommand{\U}[1]{\protect\rule{.1in}{.1in}}
\IEEEoverridecommandlockouts
\newtheorem{innercustomthm}{Theorem}
\newenvironment{customthm}[1]
  {\renewcommand\theinnercustomthm{#1}\innercustomthm}
  {\endinnercustomthm}
\newtheorem{innercustomlem}{Lemma}
\newenvironment{customlem}[1]
  {\renewcommand\theinnercustomlem{#1}\innercustomlem}
  {\endinnercustomlem}

\newtheorem{innercustomass}{Assumption}
\newenvironment{customass}[1]
  {\renewcommand\theinnercustomass{#1}\innercustomass}
  {\endinnercustomass}

\newtheorem{innercustompro}{Proposition}
\newenvironment{custompro}[1]
  {\renewcommand\theinnercustompro{#1}\innercustompro}
  {\endinnercustompro}
\allowdisplaybreaks[0]

\usepackage[normalem]{ulem}
\usepackage{color}
\usepackage{bbm}
\begin{document}

\title{{\LARGE \textbf{Optimal Control of Connected Automated Vehicles at Urban
Traffic Intersections: A Feasibility Enforcement Analysis }}}
\author{Yue Zhang, Christos G. Cassandras, Andreas A. Malikopoulos 
\thanks{This research was supported by US
Department of Energy's SMART Mobility Initiative. The work of Cassandras and
Zhang is supported in part by NSF under grants CNS- 1239021, ECCS-1509084, and
IIP-1430145, by AFOSR under grant FA9550-15-1-0471, and by a grant from the
MathWorks.} \thanks{Y. Zhang and C.G. Cassandras are with the Division of
Systems Engineering and Center for Information and Systems Engineering, Boston
University, Boston, MA 02215 USA (e-mail: joycez@bu.edu; cgc@bu.edu).}
\thanks{A.A. Malikopoulos is with the Department of Mechanical Engineering, University of Delaware, Newark, DE 19716 USA (email: andreas@udel.edu).} }
\maketitle

\begin{abstract}
Earlier work has established a decentralized optimal control framework for
coordinating on line a continuous flow of Connected Automated Vehicles (CAVs)
entering a \textquotedblleft control zone\textquotedblright\ and crossing two
adjacent intersections in an urban area. A solution, when it exists, allows
the vehicles to minimize their fuel consumption while crossing the
intersections without the use of traffic lights, without creating congestion,
and under the hard safety constraint of collision avoidance. We establish the
conditions under which such solutions exist and show that they can be enforced
through an appropriately designed \textquotedblleft feasibility enforcement
zone\textquotedblright\ that precedes the control zone. The proposed solution
and overall control architecture are illustrated through simulation.

\end{abstract}

\thispagestyle{empty} \pagestyle{empty}



\section{Introduction}

Connected and automated vehicles (CAVs) provide significant new opportunities
for improving transportation safety and efficiency using inter-vehicle as well
as vehicle-to-infrastructure communication \cite{Li2014}. To date, traffic
lights are the prevailing method used to control the traffic flow through an
intersection. More recently, however, data-driven approaches have been
developed leading to online adaptive traffic light control as in
\cite{Fleck2015}. Aside from the obvious infrastructure cost and the need for
dynamically controlling green/red cycles, traffic light systems also lead to
problems such as significantly increasing the number of rear-end collisions at
an intersection. These issues have provided the motivation for drastically new
approaches capable of providing a smoother traffic flow and more
fuel-efficient driving while also improving safety.

The advent of CAVs provides the opportunity for such new approaches. Dresner
and Stone \cite{Dresner2004} proposed a scheme for automated vehicle
intersection control based on the use of reservations whereby a centralized
controller coordinates a crossing schedule based on requests and information
received from the vehicles located inside some communication range. The main
challenges in this case involve possible deadlocks and heavy communication
requirements which can become critical. There have been numerous other efforts
reported in the literature based on such a reservation scheme
\cite{Dresner2008,DeLaFortelle2010,Huang2012}.

Increasing the throughput of an intersection is one desired goal which can be
achieved through the travel time optimization of all vehicles located within a
radius from the intersection. Several efforts have focused on minimizing
vehicle travel time under collision-avoidance constraints
\cite{Li2006,Yan2009,Zohdy2012, Zhu2015}. Lee and Park \cite{Lee2012} proposed
a different approach based on minimizing the overlap in the position of
vehicles inside the intersection rather than their arrival times. Miculescu
and Karaman \cite{Miculescu2014a} used queueing theory and modeled an
intersection as a polling system where vehicles are coordinated to cross
without collisions. There have been also several research efforts to address
the problem of vehicle coordination at intersections within a decentralized
control framework. A detailed discussion of the research in this area reported
in the literature to date can be found in \cite{Rios-Torres}.

Our earlier work \cite{ZhangMalikopoulosCassandras2016} has established a
decentralized optimal control framework for coordinating online a continuous
flow of CAVs crossing two adjacent intersections in an urban area. We refer to
an approach as \textquotedblleft centralized\textquotedblright\ if there is at
least one task in the system that is globally decided for all vehicles by a
single central controller. In contrast, in a \textquotedblleft
decentralized\textquotedblright\ approach, a \emph{coordinator}\ may be used
to handle or distribute information available in the system without, however,
getting involved in any control task. The framework in
\cite{ZhangMalikopoulosCassandras2016} solves an optimal control problem for
each CAV entering a specified \emph{Control Zone} (CZ) which subsequently
regulates the acceleration/deceleration of the CAV. The optimal control
problem involves hard safety constraints, including rear-end collision
avoidance. These constraints make it nontrivial to ensure the existence of a
feasible solution to this problem. In fact, it is easy to check that the
rear-end collision avoidance constraints cannot be guaranteed to hold
throughout the CZ under an optimal solution unless the initial conditions
(time and speed) of each CAV entering the CZ satisfy certain conditions. It
is, therefore, of fundamental importance to determine these feasibility
conditions and ensure that they can be satisfied.

The contributions of this paper are twofold. First, we study the feasibility
conditions required to guarantee a solution of the optimal control problem for
each CAV; these are expressed in terms of a feasible region defined in the
space of the CAV's speed and arrival time at the CZ. Second, we introduce a
\emph{Feasibility Enforcement Zone} (FEZ) which precedes the CZ and within
which a CAV is controlled with the goal of attaining a point in the feasible
region determined by the current state of the CZ. This subsequently guarantees
that all required constraints are satisfied when the CAV enters the CZ under
an associated optimal control. We emphasize again that the benefits of an
optimal controller maximizing throughput and minimizing fuel consumption can
only be realized subject to \emph{ensuring feasible initial conditions} to the
optimization problem under consideration.

The structure of the paper is as follows. In Section II, we review the model
in \cite{ZhangMalikopoulosCassandras2016} and its generalization in
\cite{Malikopoulos2016}. In Section III, we present the CAV coordination
framework and associated optimal control problems and solutions considering
control/state constraints. In Section IV, we carry out the analysis necessary
to identify a feasible region for the initial conditions of each CAV when
entering the CZ. In Section V, we develop a design procedure for the FEZ and
in Section VI, we include simulation results. We offer concluding remarks in
Section VII.


\section{The Model}

\begin{figure}[ptb]
\centering
\includegraphics[width=2.4 in]{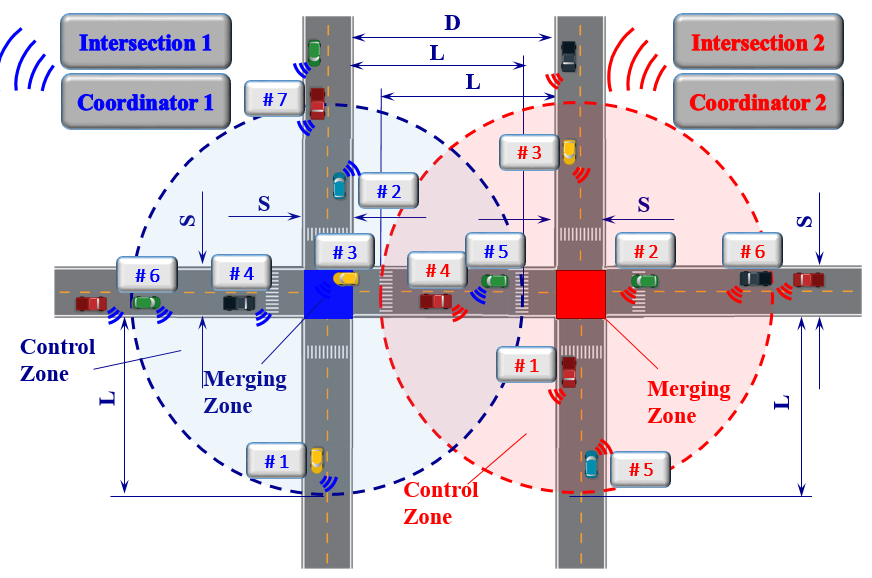} \caption{Two intersections
with connected and automated CAVs.}%
\label{fig:intersection}%
\end{figure}

We briefly review the model introduced in
\cite{ZhangMalikopoulosCassandras2016} and \cite{Malikopoulos2016} where there
are two intersections, 1 and 2, located within a distance $D$ (Fig.
\ref{fig:intersection}). The region at the center of each intersection, called
\emph{Merging Zone} (MZ), is the area of potential lateral CAV collision.
Although it is not restrictive, this is taken to be a square of side $S$. Each
intersection has a \emph{Control Zone} (CZ) and a coordinator that can
communicate with the CAVs traveling within it. The distance between the entry
of the CZ and the entry of the MZ is $L>S$, and it is
assumed to be the same for all entry points to a given CZ.

Let 
$M_{z}(t)\in\mathbb{N}$ be the cumulative number of CAVs which have entered the CZ and formed a first-in-first-out (FIFO) queue 
by time $t$, $z = 1,2$. When a CAV reaches the CZ of intersection $z$, the coordinator
assigns it an integer value $i=M_{z}(t)+1$. If two or more
CAVs enter a CZ at the same time, then the corresponding coordinator
selects randomly the first one to be assigned the value $M_{z}(t)+1$. In the region between the exit
point of a MZ and the entry point of the subsequent CZ, the CAVs cruise
with the speed they had when they exited that MZ.

For simplicity, we assume that
each CAV is governed by second order dynamics%
\begin{equation}
\dot{p}_{i}=v_{i}(t)\text{, }~p_{i}(t_{i}^{0})=0\text{; }~\dot{v}_{i}%
=u_{i}(t)\text{, }v_{i}(t_{i}^{0})\text{ given}\label{eq:model2}%
\end{equation}
where $p_{i}(t)\in\mathcal{P}_{i}$, $v_{i}(t)\in\mathcal{V}_{i}$, and
$u_{i}(t)\in\mathcal{U}_{i}$ denote the position, i.e., travel distance since the entry of the CZ, speed and
acceleration/deceleration (control input) of each CAV $i$. These dynamics are
in force over an interval $[t_{i}^{0},t_{i}^{f}]$, where $t_{i}^{0}$ is the
time that CAV $i$ enters the CZ and $t_{i}^{f}$ is the
time that it exits the MZ of intersection $z$.

To ensure that the control input and vehicle speed are within a given
admissible range, the following constraints are imposed:
\begin{equation}%
\begin{split}
u_{i,min}  &  \leqslant u_{i}(t) \leqslant u_{i,max},\quad\text{and}\\
0  &  \leqslant v_{min} \leqslant v_{i}(t) \leqslant v_{max},\quad\forall
t\in\lbrack t_{i}^{0},t_{i}^{m}],
\end{split}
\label{speed_accel constraints}%
\end{equation}
To ensure the absence of any rear-end collision throughout the CZ, we impose
the \emph{rear-end safety} constraint
\begin{equation}
s_{i}(t)=p_{k}(t)-p_{i}(t) \geqslant\delta,\quad\forall t\in\lbrack t_{i}%
^{0},t_{i}^{m}] \label{rearend}%
\end{equation}
where $\delta$ is the \emph{minimal safe distance} allowable and $k$ is the CAV physically ahead of $i$.

As part of safety considerations, we impose the following assumption (which
may be relaxed if necessary):


\begin{customass}{1}
The speed of the CAVs inside the MZ is constant, i.e., $v_{i}(t)=v_{i}%
(t_{i}^{m})=v_{i}(t_{i}^{f})$, $\ \forall t\in\lbrack t_{i}^{m},t_{i}^{f}]$,
where $t_{i}^{m}$ is the time that CAV $i$ enters the MZ
of the intersection. \label{ass:4} This implies that
\end{customass}

\begin{equation}
t_{i}^{f}=t_{i}^{m}+\frac{S}{v_{i}(t_{i}^{m})}. \label{eq:time}
\end{equation}

The objective of each CAV is to derive an optimal acceleration/deceleration in
terms of fuel consumption over the time interval $[t_{i}^{0},t_{i}^{m}]$ while
avoiding congestion between the two intersections. In addition, we impose hard
constraints so as to avoid either rear-end collision, or lateral collision inside
the MZ. In fact, it is shown in \cite{Malikopoulos2016} that the centralized
throughput maximization problem is equivalent to a set of decentralized
problems whereby each CAV minimizes its fuel consumption as long as the safety
constraints applying to it are satisfied. Thus, in what follows, we focus on
these decentralized problems and their associated safety constraints.


\section{Vehicle Coordination and Control}


\subsection{Decentralized Control Problem Formulation}

\label{sec:3a}

Since the coordinator is not involved in any decision on the vehicle control, we can
formulate $M_{1}(t)$ and $M_{2}(t)$ decentralized tractable problems for
intersection 1 and 2 respectively that may be solved on line. 
When a CAV enters a CZ, $z=1,2$, it is assigned a pair $(i,j)$ from the coordinator, where $i=M_{z}(t)+1$ is a unique index and $j$
indicates the positional relationship between CAVs $i-1$ and $i$. As formally defined in
\cite{ZhangMalikopoulosCassandras2016}, with respect to CAV $i$, CAV $i-1$ belongs to one and only one of the four following subsets:
 $(i)$ $\mathcal{R}%
_{i}^{z}(t)$ contains all CAVs traveling on the same road as $i$ and towards
the same direction but on different lanes, $(ii)$ $\mathcal{L}_{i}^{z}(t)$
contains all CAVs traveling on the same road and lane as CAV $i$, $(iii)$
$\mathcal{C}_{i}^{z}(t)$ contains all CAVs traveling on different roads from
$i$ and having destinations that can cause lateral collision at the MZ, and $(iv)$ $\mathcal{O}_{i}^{z}(t)$ contains all CAVs traveling on
the same road as $i$ and opposite destinations that cannot, however, cause
collision at the MZ. Note that the FIFO structure of this queue
implies the following condition:
\begin{equation}
t_{i}^{m}\geqslant t_{i-1}^{m},~i>1.
\label{eq:fifo}%
\end{equation}

Under the
assumption that each CAV $i$ has proximity sensors and can observe and/or
estimate local information that can be shared with other CAVs, we define its
\emph{information set} $Y_{i}(t)$, $t\in\lbrack t_{i}^{0},t_{i}^{f}]$, as
\begin{equation}
Y_{i}(t)\triangleq\Big\{p_{i}(t),v_{i}(t),\mathcal{Q}_{j}^{z},j=1,\ldots
,4,z=1,2,s_{i}(t),t_{i}^{m}\Big\},\label{InfoSet}%
\end{equation}
where $p_{i}(t),v_{i}(t)$ are the position and speed of CAV $i$ inside the
\textit{CZ }it belongs to, and $\mathcal{Q}_{j}^{z}\in\{\mathcal{R}_{i}%
^{z}(t),$ $\mathcal{L}_{i}^{z}(t),$ $\mathcal{C}_{i}^{z}(t),$ $\mathcal{O}%
_{i}^{z}(t)\},$ $z=1,2,$ is the subset assigned to CAV $i$ by the
coordinator. The fourth element in $Y_{i}(t)$ is $s_{i}(t)=p_{k}(t)-p_{i}(t)$,
the distance between CAV $i$ and CAV $k$ which is immediately
ahead of $i$ in the same lane (the index $k$ is made available to $i$ by the
coordinator). The last element above, $t_{i}^{m}$, is the time targeted for
CAV $i$ to enter the MZ, whose evaluation is discussed next. Note that
once CAV $i$ enters the CZ, then all information in $Y_{i}(t)$ becomes
available to $i$.

The time $t_{i}^{m}$ that CAV $i$ is required to enter the MZ is based on
maximizing the intersection throughput while satisfying (\ref{eq:fifo}) and
the constraints for avoiding rear-end and lateral collision in the MZ. There
are three cases to consider regarding $t_{i}^{m}$, depending on the value of
$\mathcal{Q}_{j}^{z}$:

\emph{Case 1}: $(i-1)\in\mathcal{R}_{i}^{z}(t)\cup\mathcal{O}_{i}^{z}(t):$ in this
case, none of the safety constraints can become active while $i$ and $i-1$ are
in the CZ or MZ. This allows CAV $i$ to minimize its time in the CZ while
preserving the FIFO queue through $t_{i}^{m} \geqslant t_{i-1}^{m}, 
~ i>1$. Therefore, it is obvious that we should set
\begin{equation}
t_{i}^{m}=t_{i-1}^{m} \label{eq:lem1a}%
\end{equation}
and since CAV speeds inside the MZ are constant (Assumption \ref{ass:4}), both
$i-1$ and $i$ will also be exiting the MZ at the same time by setting%
\begin{equation}
v_{i}^m =v_{i-1}^m. \label{Case1_speed}%
\end{equation}
where $v_i^m = v_i(t_i^m)$ and $v_{i-1}^m = v_{i-1}(t_{i-1}^m)$. Note that, by Assumption \ref{ass:4}, $v_{i}(t)=v_{i-1}(t)$ for all
$t\in\lbrack t_{i}^{m},t_{i}^{f}]$.

\emph{Case 2}: $(i-1)\in\mathcal{L}_{i}^{z}(t):$ in this case, only the rear-end
collision constraint \eqref{rearend} can become active. In order to minimize
the time CAV $i$ spends in the CZ by ensuring that \eqref{rearend} is
satisfied over $t\in\lbrack t_{i}^{m},t_{i-1}^{f}]$ while $v_{i-1}(t)$ is
constant (Assumption \ref{ass:4}), we set%
\begin{equation}
t_{i}^{m}=t_{i-1}^{m}+\frac{\delta}{v_{i-1}^{m}}, \label{eq:lem1b}%
\end{equation}
and $v_i^m$ as in \eqref{Case1_speed}.

\emph{Case 3}: $(i-1)\in\mathcal{C}_{i}^{z}(t):$ in this case, only the lateral
collision may occur. Hence, CAV $i$ is allowed to enter the MZ only when
CAV $i-1$ exits from it. To minimize the time CAV $i$ spends in the CZ
while ensuring that the lateral collision avoidance is satisfied over
$t\in\lbrack t_{i}^{m},t_{i-1}^{f}]$, we set%
\begin{equation}
t_{i}^{m}=t_{i-1}^{m}+\frac{S}{v_{i-1}^{m}}, \label{eq:lem1c}%
\end{equation}
and $v_i^m$ as in \eqref{Case1_speed}.

It follows from (\ref{eq:lem1a}) through (\ref{eq:lem1c}) that $t_{i}^{m}$ is
always recursively determined from $t_{i-1}^{m}$ and $v_{i-1}^m$.
Similarly, $v_i^m$ depends only on $v_{i-1}^{m}$.

Although (\ref{eq:lem1a}), (\ref{eq:lem1b}), and (\ref{eq:lem1c}) provide a
simple recursive structure for determining $t_{i}^{m}$, the presence of the
control and state constraints \eqref{speed_accel constraints} may prevent
these values from being admissible. This may happen by
\eqref{speed_accel constraints} becoming active at some internal point during
an optimal trajectory (see \cite{Malikopoulos2016} for details). In addition,
however, there is a global lower bound to $t_{i}^{f}$, hence also $t_{i}^{m}$
through (\ref{eq:time}), which depends on $t_{i}^{0}$ and on whether CAV $i$
can reach $v_{max}$ prior to $t_{i-1}^{m}$ or not: $(i)$ If CAV $i$ enters the
CZ at $t_{i}^{0}$, accelerates with $u_{i,max}$ until it reaches $v_{max}$ and
then cruises at this speed until it leaves the MZ at time $t_{i}^{1}$, it was
shown in \cite{ZhangMalikopoulosCassandras2016} that
\begin{equation}
t_{i}^{1}=t_{i}^{0}+\frac{L+S}{v_{max}}+\frac{(v_{max}-v_{i}^{0})^{2}%
}{2u_{i,max}v_{max}}.
\end{equation}
$(ii)$ If CAV $i$ accelerates with $u_{i,max}$ but reaches the MZ at $t_{i}^{m}$
with speed $v_{i}^{m}<v_{max}$, it was shown in
\cite{ZhangMalikopoulosCassandras2016} that
\begin{equation}
t_{i}^{2}=t_{i}^{0}+\frac{v_{i}(t_{i}^{m})-v_{i}^{0}}{u_{i,max}}+\frac{S}%
{v_{i}^{m}},
\end{equation}
where $v_{i}(t_{i}^{m})=\sqrt{2Lu_{i,max}+(v_{i}^{0})^{2}}$. Thus,
\begin{equation}
t_{i}^{c}=t_{i}^{1} \mathbbm{1}_{ v_i^m = v_{max}} + t_{i}^{2} (1 - \mathbbm{1}_{ v_i^m = v_{max}} ) \nonumber
\end{equation}
is a lower bound of $t_{i}^{f}$ regardless
of the solution of the problem. Therefore, we can summarize the recursive
construction of $t_{i}^{f}$ over $i=1,\ldots,M_z(t)$ as follows:%
\begin{equation}
t_{i}^{f}=\left\{
\begin{array}
[c]{ll}%
t_{1}^{f}, & \mbox{if $i=1$},\\
\text{max }\{t_{i-1}^{f},t_{i}^{c}\}, & \text{if }i-1\in\mathcal{R}_{i}^{z}%
(t)\cup\mathcal{O}_{i}^{z}(t)\\
\text{max }\{t_{i-1}^{f}+\frac{\delta}{v_{i}(t_{i-1}^{f})},t_{i}^{c}\}, &
\mbox{if $i-1\in\mathcal{L}_{i}^{z}$},\\
\text{max }\{t_{i-1}^{f}+\frac{S}{v_{i}(t_{i-1}^{f})},t_{i}^{c}\}, &
\mbox{if $i-1\in\mathcal{C}_{i}^{z}$},
\end{array}
\right.  \label{def:tf}%
\end{equation}
where $t_{i}^{m}$ can be evaluated from $t_{i}^{f}$ through (\ref{eq:time}),
and thus, it is always feasible.


Note that at each time $t$, each CAV $i$ communicates with the preceding CAV $i-1$ in the queue and accesses the
values of $t_{i-1}^{f}$, $v_{i-1}(t_{i-1}^{f})$, $\mathcal{Q}_{j}^{z}$,
$j=1,\dots,4$, $z=1,2$ from its information set in (\ref{InfoSet}). This is
necessary for $i$ to compute $t_{i}^{f}$ appropriately and satisfy
\eqref{def:tf} and (\ref{rearend}). The following result is established in
\cite{ZhangMalikopoulosCassandras2016} to formally assert the iterative
structure of the sequence of decentralized optimal control problems:

\begin{customlem}{1}
The decentralized communication structure aims for each CAV $i$ to solve an
optimal control problem for $t\in\lbrack t_{i}^{0},t_{i}^{m}]$ the solution of
which depends only on the solution of CAV $i$-1.
\end{customlem}

The decentralized optimal control problem for each CAV approaching either
intersection is formulated so as to minimize the $L^{2}$-norm of its control
input (acceleration/deceleration). It has been shown in \cite{Rios-Torres2015}
that there is a monotonic relationship between fuel consumption for each
CAV $i$, and its control input $u_{i}$. Therefore, we formulate the
following problem for each $i$:
\begin{gather}
\min_{u_{i}\in U_{i}}\frac{1}{2}\int_{t_{i}^{0}}^{t_{i}^{m}}K_{i}\cdot
u_{i}^{2}~dt\nonumber\\
\text{subject to}:\eqref{eq:model2},(\ref{speed_accel constraints}%
),(\ref{eq:time}),\eqref{def:tf},\text{ }p_{i}(t_{i}^{0})=0\text{, }%
p_{i}(t_{i}^{m})=L,\label{eq:decentral}\\
z=1,2, \text{and given }t_{i}^{0}\text{, }v_{i}(t_{i}^{0}).\nonumber
\end{gather}
where $K_{i}$ is a factor to capture CAV diversity (for simplicity we set
$K_{i}=1$ for the rest of this paper). Note that this formulation does not include the safety constraint (\ref{rearend}).


\subsection{Analytical solution of the decentralized optimal control problem}

\label{sec:4b}

An analytical solution of problem \eqref{eq:decentral} may be obtained through
a Hamiltonian analysis. The presence of constraints (\ref{speed_accel constraints}) and
\eqref{def:tf} complicates this analysis. Assuming that all constraints are
satisfied upon entering the CZ and that they remain inactive throughout
$[t_{i}^{0},t_{i}^{m}]$, a complete solution was derived in
\cite{Rios-Torres2015} and \cite{Rios-Torres2} for highway on-ramps, and in
\cite{ZhangMalikopoulosCassandras2016} for two adjacent intersections. This
solution is summarized next (the complete solution including any constraint
(\ref{speed_accel constraints}) becoming active is given in
\cite{Malikopoulos2016}). The optimal control input
(acceleration/deceleration) over $t\in\lbrack t_{i}^{0},t_{i}^{m}]$ is given
by
\begin{equation}
u_{i}^{\ast}(t)=a_{i}t+b_{i}\label{eq:20}%
\end{equation}
where $a_{i}$ and $b_{i}$ are constants. Using (\ref{eq:20}) in the CAV
dynamics \eqref{eq:model2} we also obtain the optimal speed and position:
\begin{equation}
v_{i}^{\ast}(t)=\frac{1}{2}a_{i}t^{2}+b_{i}t+c_{i}\label{eq:21}%
\end{equation}%
\begin{equation}
p_{i}^{\ast}(t)=\frac{1}{6}a_{i}t^{3}+\frac{1}{2}b_{i}t^{2}+c_{i}%
t+d_{i},\label{eq:22}%
\end{equation}
where $c_{i}$ and $d_{i}$ are constants of integration. The constants $a_{i}$,
$b_{i}$, $c_{i}$, $d_{i}$ can be computed by using the given initial and final
conditions.
The interdependence of the two intersections, i.e., the coordination of
CAVs at the MZ of one intersection which affects the behavior of CAV
coordination of the other MZ, is discussed in
\cite{ZhangMalikopoulosCassandras2016}.

We note that the control of CAV $i$ actually remains unchanged until an
\textquotedblleft event\textquotedblright\ occurs that affects its behavior.
Therefore, the time-driven controller above can be replaced by an event-driven
one without affecting its optimality properties under conditions described in
\cite{Zhong2010}.


As already mentioned, the analytical solution (\ref{eq:20}) is only
valid as long as all initial conditions satisfy (\ref{speed_accel constraints}) and
\eqref{def:tf} and these constraints continue to be satisfied throughout
$[t_{i}^{0},t_{i}^{m}]$. Otherwise, the solution needs to be modified as described in \cite{Malikopoulos2016}.

Recall that the constraint (\ref{rearend}) is not
included in (\ref{eq:decentral}) and it is a much more challenging matter. To
deal with this, we proceed as follows. First, we analyze under what initial
conditions $(t_{i}^{0},v_{i}^{0})$ the constraint is violated upon CAV $i$
entering the CZ. This defines a \emph{feasibility region} in the $(t_{i}%
^{0},v_{i}^{0})$ space which we denote by $\mathcal{F}_{i}$. Assuming the CAV
has initial conditions which are feasible, we then derive a condition under
which the CAV's state maintains feasibility over $[t_{i}^{0},t_{i}^{m}]$.
Finally, we explore how to enforce feasibility at the time of CZ entry, i.e.,
enforcing the condition $(t_{i}^{0},v_{i}^{0})\in\mathcal{F}_{i}$. This is
accomplished by introducing a \emph{Feasibility Enforcement Zone} (FEZ) which
precedes the CZ. If the FEZ is properly designed, we show that $(t_{i}%
^{0},v_{i}^{0})\in\mathcal{F}_{i}$ can be ensured.

\section{Feasibility enforcement analysis}

\begin{figure}[ptb]
\centering
\includegraphics[width=2 in]{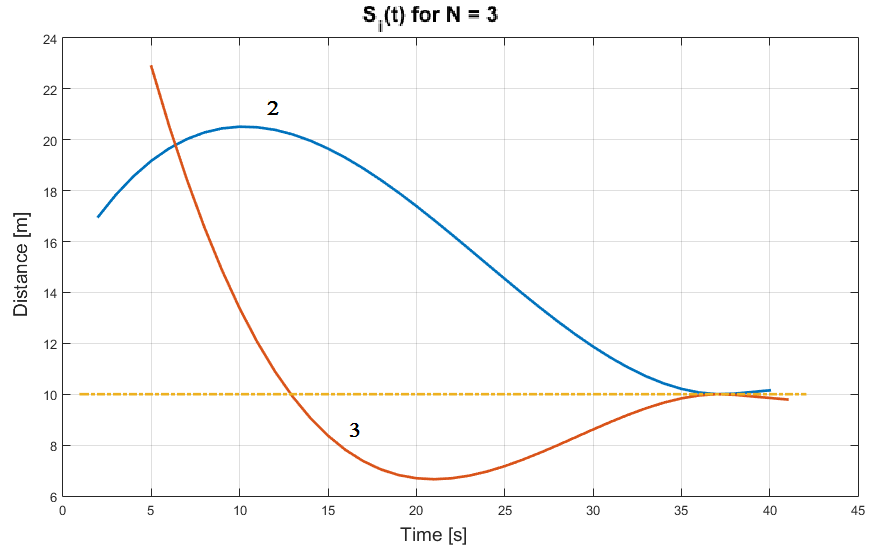}\caption{Example of safety
constraint violation by CAV 3 when $\delta$ = 10.}%
\label{viocase}%
\end{figure}

We begin with a simple example of how the safety constraint (\ref{rearend})
may be violated under the optimal control (\ref{eq:20}). This is illustrated
in Fig. \ref{viocase} with $\delta=10$ for two CAVs that follow each other
into the same lane in the CZ. We can see that while (\ref{rearend}) is
eventually satisfied over the MZ, due to the constraints imposed on the solution of
\eqref{eq:decentral} through (\ref{def:tf}), the controller (\ref{eq:20}) is
unable to maintain (\ref{rearend}) throughout the CZ. What is noteworthy in
Fig. \ref{viocase} is that (\ref{rearend}) is violated by CAV 3 at an interval
which is \emph{interior} to $[t_{3}^{0},t_{3}^{m}]$, i.e., the form of the
optimal control solution (\ref{eq:20}) causes this violation even though the
constraint is initially satisfied at $t_{3}^{0}=5$ in Fig. \ref{viocase}.

Recall that we use $k$ to denote the CAV physically preceding $i$ on the same lane in the
CAV, and $i-1$ is the CAV preceding $i$ in the FIFO queue associated with the
CAV, we have the following theorem.

\begin{customthm}{1}
There exists a nonempty feasible region $\mathcal{F}_{i}\subset R^{2}$ of
initial conditions $(t_{i}^{0},v_{i}^{0})$ for CAV $i$ such that, under the
decentralized optimal control, $s_{i}(t)\geqslant\delta$ for all $t\in\lbrack
t_{i}^{0},t_{i}^{m}]$ given the initial and final conditions $t_{k}^{0}%
,v_{k}^{0},t_{k}^{m},v_{k}^{m}$ of CAV $k$. \label{theo:3}
\end{customthm}

 
\emph{Proof:} To prove the existence of the feasible region, there are two
cases to consider, depending on whether any state or control constraint for
either CAV $i$ or $k$ becomes active in the CZ.

\emph{Case 1:} No state or control constraint is active for either $k$ or $i$
over $[t_{i}^{0},t_{i}^{m}]$. By using \eqref{eq:21}, \eqref{eq:22} at $t$ and
$t_{i}^{m}$, and the definition $s_{i}(t)=p_{k}(t)-p_{i}(t)$, under optimal
control we can write%
\begin{gather}
s_{i}(t; t_{i}^{0}, v_{i}^{0}) = s_{i}(t, t_{i}^{m},v_{i}^{m},t_{k}^{0}%
,v_{k}^{0},t_{k}^{m},v_{k}^{m}; t_{i}^{0}, v_{i}^{0})\nonumber\\
= A(t, t_{i}^{m},v_{i}^{m},t_{k}^{0},v_{k}^{0},t_{k}^{m},v_{k}^{m}; t_{i}^{0},
v_{i}^{0})t^{3}\nonumber\\
+ B(t, t_{i}^{m},v_{i}^{m},t_{k}^{0},v_{k}^{0},t_{k}^{m},v_{k}^{m}; t_{i}^{0},
v_{i}^{0})t^{2}\nonumber\\
+ C(t, t_{i}^{m},v_{i}^{m},t_{k}^{0},v_{k}^{0},t_{k}^{m},v_{k}^{m}; t_{i}^{0},
v_{i}^{0})t\nonumber\\
+ D(t, t_{i}^{m},v_{i}^{m},t_{k}^{0},v_{k}^{0},t_{k}^{m},v_{k}^{m}; t_{i}^{0},
v_{i}^{0}), \label{si(t)}%
\end{gather}
where $A$, $B$, $C$ and $D$ are functions defined over $t\in\lbrack t_{i}%
^{0},t_{i}^{m}]$. Recall that CAV $k$ is cruising in the MZ, so that
(\ref{eq:20}) through (\ref{eq:22}) do not apply for $k$ over $[t_{k}%
^{m},t_{i}^{m}]$ leading to different expressions for $A$, $B$, $C$ and $D$.
Therefore, we consider two further subcases, one for $[t_{i}^{0},t_{k}^{m}]$
and the other for $[t_{k}^{m},t_{i}^{m}]$. For ease of notation, in the sequel
we replace $(t_{i}^{0},v_{i}^{0})$ by $(\tau,\upsilon)$.

\emph{Case 1.1:} $t\in\lbrack t_{i}^{0},t_{k}^{m}]$. In this case,
$s_{i}(t;\tau,\upsilon)$ is a cubic polynomial inheriting the cubic structure
of (\ref{eq:22}). We can solve \eqref{eq:21}, \eqref{eq:22} for the
cofficients $a_{k}$, $b_{k}$, $c_{k}$, $d_{k}$, $a_{i}$, $b_{i}$, $c_{i}$ and
$d_{i}$ using the initial and final conditions of CAVs $k$ and $i$. Then,
denoting $A$, $B$, $C$ and $D$ as $A_{1}(\tau,v)$, $B_{1}(\tau,v)$,
$C_{1}(\tau,v)$ and $D_{1}(\tau,v)$ for $t\in\lbrack t_{i}^{0},t_{k}^{m}]$,
these are explicitly given by%
\begin{equation}
\begin{aligned}
A_{1}(\tau,\upsilon) & = \frac{1}{(t_{k}^{0}-t_{k}^{m})^{3}}(2L + (v_{k}%
^{m}+v_{k}^{0})(t_{k}^{0}-t_{k}^{m}))\nonumber\\
& - \frac{1}{(\tau-t_{i}^{m})^{3}}(2L + (v_{i}^{m}+\upsilon)(\tau-t_{i}%
^{m})),\nonumber \\
B_{1}(\tau,\upsilon) & = - \frac{1}{(t_{k}^{0}-t_{k}^{m})^{3}}[3L(t_{k}%
^{0}+t_{k}^{m})\nonumber\\
&+(v_{k}^{0}(t_{k}^{0}+2t_{k}^{m}) +v_{k}^{m}(2t_{k}^{0}+t_{k}^{m}))(t_{k}%
^{0}-t_{k}^{m})]\nonumber\\
&+ \frac{1}{(\tau-t_{i}^{m})^{3}} [3L(\tau+t_{i}^{m})\nonumber\\
&+ (\upsilon(\tau+2t_{i}^{m})+v_{i}^{m}(2\tau+t_{i}^{m}))(\tau-t_{i}%
^{m})],\nonumber \\
C_{1}(\tau,\upsilon) &= \frac{1}{(t_{k}^{0}-t_{k}^{m})^{3}}[6t_{k}^{0}
t_{k}^{m} L + [(v_{k}^{0}((t_{k}^{m})^{2}+2t_{k}^{0} t_{k}^{m})\nonumber\\
& + v_{k}^{m}((t_{k}^{0})^{2}+2t_{k}^{m} t_{k}^{0}))](t_{k}^{0}-t_{k}%
^{m})]\nonumber\\
&- \frac{1}{(\tau-t_{i}^{m})^{3}}[6\tau t_{i}^{m} L + [(\upsilon((t_{i}%
^{m})^{2}+2\tau t_{i}^{m})\nonumber\\
&+ v_{i}^{m}((\tau)^{2}+2t_{i}^{m} \tau))](\tau-t_{i}^{m})],\nonumber \\
D_{1}(\tau,\upsilon) &= \frac{1}{(t_{k}^{0}-t_{k}^{m})^{3}} [L((t_{k}^{0}%
)^{3}-3(t_{k}^{0})^{2} t_{k}^{m})\nonumber\\
&- ( v_{k}^{0} t_{k}^{0}(t_{k}^{m})^{2}+ v_{k}^{m} (t_{k}^{0})^{2} t_{k}%
^{m})(t_{k}^{0}-t_{k}^{m})]\nonumber\\
&- \frac{1}{(\tau-t_{i}^{m})^{3}} [ L((\tau)^{3}-3(\tau)^{2} t_{i}%
^{m})\nonumber\\
&- (\upsilon\tau(t_{i}^{m})^{2}+v_{i}^{m} (\tau)^{2} t_{i}^{m})(\tau-t_{i}%
^{m})].\label{ABCDcoff}%
\end{aligned}
\end{equation}

Aside from $(\tau, v)$, all remaining arguments are known to CAV $i$ and can
be determined. Hence, $s_{i}(t;\tau,v)$ varies only with $t$ and $(\tau,v)$.
First, observing that the first half of each of the coefficient expressions in
\eqref{ABCDcoff} (which is derived by solving \eqref{eq:21} and \eqref{eq:22}
for CAV $k$) is a constant fully determined by information provided by CAV
$k$, we can rewrite these as $K_{A_{1}}$, $K_{B_{1}}$, $K_{C_{1}}$, $K_{D_{1}%
}$. Therefore, $p_{k}^{\ast}(t)$ in \eqref{eq:22} can be expressed as
\begin{equation}
p_{k}^{\ast}(t)=K_{A_{1}}t^{3}+K_{B_{1}}t^{2}+K_{C_{1}}t+K_{D_{1}}.
\label{pkt}%
\end{equation}
Next, the second half of the coefficients can be expressed through polynomials
in either $\tau$ or $\upsilon$ explicitly derived by solving \eqref{eq:21} and
\eqref{eq:22}  for CAV $i$. We will use the notation $P_{X,n}(\tau)$,
$P_{X,n}(\upsilon)$ to represent polynomials of degree $n=1,2,3$ and
$X\in\{A_{1},B_{1},C_{1},D_{1}\}$. Similarly, we set $Q_{3}(\tau)=(\tau
-t_{i}^{m})^{3}$. Thus, for the coefficients in Eq. \eqref{ABCDcoff}, we get%
\begin{equation}
\begin{aligned} A_1(\tau, \upsilon) &= K_{A_1} + \frac{ P_{A_1,1}(\tau) P_{A_1,1}(\upsilon)}{Q_3(\tau)}, \\ B_1(\tau, \upsilon) &= K_{B_1} + \frac{ P_{B_1,2}(\tau)P_{B_1,1}(\upsilon)}{Q_3(\tau)}, \\ C_1(\tau, \upsilon) & = K_{C_1} + \frac{ P_{C_1,3}(\tau) + P_{C_1,2}(\tau)P_{C_1,1}(\upsilon)}{Q_3(\tau)},\\ D_1(\tau, \upsilon) & = K_{D_1} + \frac{ P_{D_1,3}(\tau) + P_{D_1,2}(\tau)P_{D_1,1}(\upsilon)}{Q_3(\tau)}. \end{aligned} \label{ABCD}%
\end{equation}
Note that $p_{k}^{\ast}(t)$ in (\ref{pkt}) involves only the $K$ terms, while
the analogous cubic polynomial for $p_{i}^{\ast}(t)$ involves only the $P$ and
$Q$ terms.

Our goal is to ensure that $s_{i}(t;\tau,\upsilon)\geqslant\delta$ for all
$t\in\lbrack\tau,t_{k}^{m}]$ (recall that $t_{i}^{0}\equiv\tau$). We can
guarantee this by ensuring that $s_{i}^{\ast}(\tau,\upsilon)\equiv\min
_{t\in\lbrack\tau,t_{k}^{m}]}\{s_{i}(t;\tau,\upsilon)\}\geqslant\delta$. Thus,
we shift our attention to the determination of $s_{i}^{\ast}(\tau,\upsilon)$.
We can obtain expressions for the first and the second derivative of
$s_{i}(t;\tau,\upsilon)$, $\dot{s}_{i}(t;\tau,\upsilon)$ and $\ddot{s}%
_{i}(t;\tau,\upsilon)$ respectively, from \eqref{si(t)}, as follows:%
\begin{equation}
\begin{aligned} \dot s_i(t;\tau,\upsilon) & = v_k(t) - v_i(t) \\ & = 3A_1(\tau, \upsilon) t^2 + 2B_1(\tau, \upsilon) t + C_1(\tau, \upsilon), \end{aligned} \label{si_dot}%
\end{equation}%
\begin{equation}
\begin{aligned} \ddot s_i(t;\tau,\upsilon) & = u_k(t) - u_i(t) \\ & = 6A_1(\tau, \upsilon) t + 2B_1(\tau, \upsilon). \end{aligned} \label{si_dotdot}%
\end{equation}
Clearly, we can determine $t_{i}^{\ast}\equiv\arg\min_{t\in\lbrack\tau
,t_{k}^{m}]}\{s_{i}(t;\tau,\upsilon)\}$ as the solution of $\dot{s}_{i}%
(t;\tau,\upsilon)=0$ with $\ddot{s}_{i}(t;\tau,\upsilon)\geqslant0$, unless
$s_{i}^{\ast}(\tau,\upsilon)$ occurs at the boundaries, i.e., $t_{i}^{\ast
}=\tau$ or $t_{i}^{\ast}=t_{k}^{m}$. Thus, there are three cases to consider:

\emph{Case 1.1.A}: $t_{i}^{\ast}=\tau$. In this case,%
\begin{gather}
s_{i}^{\ast}(\tau,\upsilon)=s_{i}(\tau;\tau,\upsilon)\label{casei}\\
=A_{1}(\tau,\upsilon)\tau^{3}+B_{1}(\tau,\upsilon)\tau^{2}+C_{1}(\tau
,\upsilon)\tau+D_{1}(\tau,\upsilon)\geqslant\delta\nonumber
\end{gather}
and we can satisfy $s_{i}(\tau,\upsilon)\geqslant\delta$ for any $\upsilon$ as
long as a feasible $\tau$ is determined. Since at $t=\tau$, we have
$p_{i}(\tau)=0$ and using the definition of $s_{i}(t)=p_{k}(t)-p_{i}(t)$ and
\eqref{pkt},
\[
s_{i}(\tau)=p_{k}^{\ast}(\tau)=K_{A_{1}}\tau^{3}+K_{B_{1}}\tau^{2}+K_{C_{1}%
}\tau+K_{D_{1}}.
\]
Observe that if $p_{k}(\tau)\geqslant\delta$, then CAV $i$ enters the CZ at a
safe distance from its preceding CAV $k$ and since $t_{i}^{\ast}=\tau$, we
have $s_{i}(t;\tau,\upsilon)\geqslant\delta$ for all $t\in\lbrack\tau
,t_{k}^{m}]$. Thus, it suffices to select
\begin{equation}
\tau\geqslant t_{k}^{\delta}\label{Case1_tau}%
\end{equation}
where $t_{k}^{\delta}$ is the smallest real root of $p_{k}(\tau)-\delta=0$.

\emph{Case 1.1.B}: $t_{i}^{\ast}=t_{k}^{m}$. In this case,%
\begin{gather}
s_{i}^{\ast}(\tau,\upsilon)=s_{i}(t_{k}^{m};\tau,\upsilon)\label{Case2_tkm}\\
=A_{1}(\tau,\upsilon)(t_{k}^{m})^{3}+B_{1}(\tau,\upsilon)(t_{k}^{m})^{2}%
+C_{1}(\tau,\upsilon)t_{k}^{m}\nonumber\\
+D_{1}(\tau,\upsilon)\geqslant\delta\nonumber
\end{gather}
Thus, the feasibility region $\mathcal{F}_{i}$ is defined by all
$(\tau,\upsilon)$ such that $s_{i}(t_{k}^{m};\tau,\upsilon)-\delta\geqslant0$
in the $(\tau,\upsilon)$ space.

\emph{Case 1.1.C}: $t_{i}^{\ast}=t_{1}\in(\tau,t_{k}^{m})$. This case only
arises if the determinant $\mathcal{D}_{i}(\tau,\upsilon)$ of (\ref{si_dot})
is positive, i.e.,
\begin{equation}
\mathcal{D}_{i}(\tau,\upsilon)=4B_{1}(\tau,\upsilon)^{2}-12A_{1}(\tau
,\upsilon)C_{1}(\tau,\upsilon)>0 \label{t1}%
\end{equation}
and we get%
\begin{equation}
t_{1}=\frac{-2B_{1}(\tau,\upsilon)\pm\sqrt{\mathcal{D}_{i}(\tau,\upsilon)}%
}{6A_{1}(\tau,\upsilon)} \label{t1root}%
\end{equation}
In addition, we must have
\begin{equation}
\tau<t_{1}<t_{k}^{m},\text{ \ \ }\dot{s}_{i}(t_{1};\tau,\upsilon)=0,\text{
\ \ }\ddot{s}_{i}(t_{1};\tau,\upsilon)\geqslant0 \label{t2}%
\end{equation}
Therefore, the feasibility region $\mathcal{F}_{i}$ is defined by all
$(\tau,\upsilon)$ such that%
\begin{equation}
\begin{aligned} s_{i}^{\ast}(\tau,\upsilon) & =s_{i}(t_{1};\tau,\upsilon)\\ =A_1(\tau,\upsilon)(t_{1})^{3}+ & B_1(\tau,\upsilon)(t_{1})^{2}+C_1(\tau,\upsilon )t_{1}\\ + &D_1(\tau,\upsilon)\geqslant\delta \end{aligned} \label{eq:11c}%
\end{equation}
in conjunction with (\ref{t1root})-(\ref{t2}).

\emph{Case 1.2:} $t\in\lbrack t_{k}^{m},t_{i}^{m}]$. Over this interval,
$v_{k}(t)=v_{k}^{m}$ by Assumption \ref{ass:4}. Therefore,
\eqref{eq:20}-\eqref{eq:22} no longer apply: \eqref{eq:20} becomes
$u_{k}^{\ast}(t)=0$, \eqref{eq:21} becomes $v_{k}^{\ast}(t)=v_{k}^{m}$ and
\eqref{eq:22} becomes $p_{k}^{\ast}(t)=L+v_{k}^{m}(t-t_{k}^{m})$. Evaluating
$s_{i}(t)=p_{k}(t)-p_{i}(t)$ in this case yields the following coefficients in
(\ref{ABCDcoff}):%

\begin{equation}
\begin{aligned} A_2(\tau,\upsilon) &= - \frac{1}{(\tau-t_i^m)^3}(2L + (v_i^m+\upsilon)(\tau-t_i^m)),\\ B_2(\tau,\upsilon) &= \frac{1}{(\tau-t_i^m)^3} [3L(\tau+t_i^m) + (\upsilon(\tau+2t_i^m) \\ & +v_i^m(2\tau+t_i^m))(\tau-t_i^m)],\\ C_2(\tau,\upsilon) & = v_k^m - \frac{1}{(\tau-t_i^m)^3}[6\tau t_i^m L + [(\upsilon((t_i^m)^2+2\tau t_i^m) \\ & + v_i^m((\tau)^2+2t_i^m \tau))](\tau-t_i^m)],\\ D_2(\tau,\upsilon) & = L - v_k^m t_k^m - \frac{1}{(\tau-t_i^m)^3} [ L((\tau)^3-3(\tau)^2 t_i^m) \\\ & - (\upsilon \tau (t_i^m)^2+v_i^m (\tau)^2 t_i^m)(\tau-t_i^m)]. \end{aligned}\label{MZphase}%
\end{equation}
It follows that $K_{A_{1}},K_{B_{1}},K_{C_{1}}$ and $K_{D_{1}}$ in
\eqref{ABCD} should be modified accordingly, giving $K_{A_{2}}=K_{B_{2}}=0$,
$K_{C_{2}}=v_{k}^{m}$ and $K_{D_{2}}=L-v_{k}^{m}t_{k}^{m}$. Since we are
assuming that no control or state constraints are active for CAV $i$, the
designated final time $t_{i}^{m}$ under optimal control satisfies
(\ref{eq:lem1b}), i.e., $s_{i}(t_{i}^{m})=\delta$. Thus, we only need to
consider the subcase where $s_{i}^{\ast}(\tau,\upsilon)$ occurs in $(t_{k}%
^{m},t_{i}^{m})$ and we have $t_{i}^{\ast}=t_{2}$, $t_{2}\in(t_{k}^{m}%
,t_{i}^{m})$. Proceeding as in \emph{Case 1.1.C}, the feasibility region
$\mathcal{F}_{i}$ is defined by all $(\tau,\upsilon)$ such that%
\begin{gather}
s_{i}^{\ast}(\tau,\upsilon)=s_{i}(t_{2};\tau,\upsilon)\\
=A_{2}(\tau,\upsilon)(t_{2})^{3}+B_{2}(\tau,\upsilon)(t_{2})^{2}+C_{2}%
(\tau,\upsilon)t_{2}\nonumber\\
+D_{2}(\tau,\upsilon)\geqslant\delta\nonumber
\end{gather}
in conjunction with (\ref{t1root})-(\ref{t2}), with $A_{1}$, $B_{1}$, $C_{1}$
and $D_{1}$ replaced by $A_{2}$, $B_{2}$, $C_{2}$ and $D_{2}$, and with
$\tau<t_{1}<t_{k}^{m}$ replaced by $t_{k}^{m}<t_{2}<t_{i}^{m}$.

\emph{Case 2:} At least one of the state and control constraints is active
over $[\tau,t_{i}^{m}]$. The analysis for this case is similar and is omitted
but it may be found in \cite{Malikopoulos2016}.

To complete the proof, we show that feasibility region $\mathcal{F}_{i}$ is
always nonempty. This is easily established by considering a point
$(\tau,\upsilon)$ such that $v_{min}<\upsilon<v_{max}$ and $\tau=t_{k}^{f}$:
since $p_{k}^{\ast}(t_{k}^{f})=L+S$ and $p_{i}^{\ast}(\tau)=0$, it follows
that $s_{i}(\tau)>S>\delta$. Obviously, any such $(\tau,\upsilon)$ is
feasible. $\blacksquare$

\begin{figure}[ptb]
\centering
\includegraphics[width = 3 in]{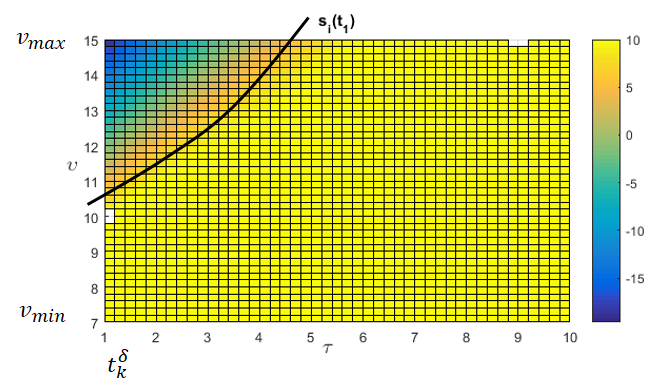} \caption{Illustration of
the feasibility region for case 1.1.C.}%
\label{case11c}%
\end{figure}

To illustrate the feasible region and provide some intuition, we give a
numerical example where \emph{Case 1.1.C} applies 
(see Fig. \ref{case11c}) with $\delta=10$, $L=400$, and CAV $k$ is the first CAV in the CZ and is
driving at the constant speed $v_{k}^{m}=10$. The colorbar in Fig.
\ref{case11c} indicates the value of $s_{i}^{\ast}(t)$ and the yellow region determined by \eqref{eq:11c}, 
represents the feasible region, while the
non-yellow region represents the infeasible region. The black curve is the
boundary between the two regions and is not linear in general. This boundary
curve shifts depending on the different cases we have considered in the proof
of Theorem \ref{theo:3}. This example also illustrates that we can always find
a nonempty feasible region since we can select points to the right of the
curve corresponding to CAV $i$ entry times in the CZ which can be arbitrarily large.

\section{Design of the feasibility enforcement zone}

Given all the information pertaining to CAVs $k$ and $i-1$, we can immediately
determine the feasible region $\mathcal{F}_{i}$ for any CAV $i$ which may
enter the CZ next. The role of the FEZ introduced prior to the CZ is to exert
a control on $i$ that ensures its initial condition $(\tau,\upsilon)$ is such
that $(\tau,\upsilon)\in\mathcal{F}_{i}$. Thus, while an optimal control is
applied to $i$ within the CZ, the control used in the FEZ is not optimal, but
it is necessary to guarantee that the subsequent optimal control is feasible.
This is similar to controllers used at gateways of communication networks in
order to \textquotedblleft smooth\textquotedblright\ incoming traffic before
applying optimal routing or scheduling policies on packets entering the
network (in our case, CAVs entering a CZ).

\begin{figure}[ptb]
\centering
\includegraphics[width = 2 in]{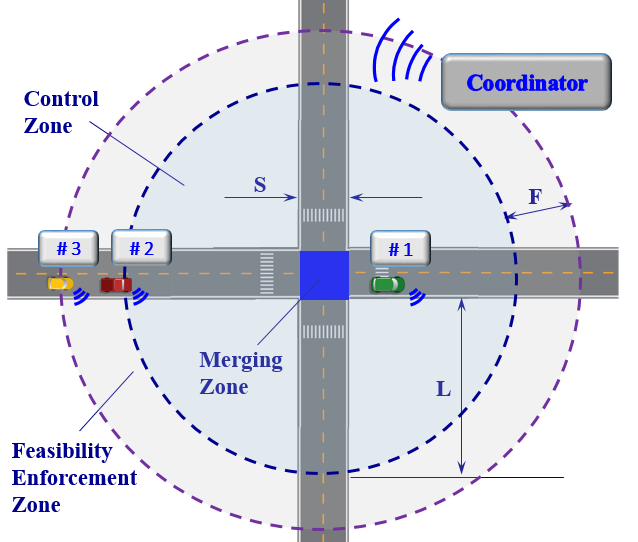} \caption{Intersection
model with feasibility enforcement zone (FEZ) added.}%
\label{fig:single}%
\end{figure}

The design of the FEZ rests on determining its length, denoted by $F_{i}$. Let
$i$ denote a CAV entering the FEZ and let $k$ denote the CAV immediately
preceding it. Let $v_{i}^{F}=v_{i}(t_{i}^{F})$ be the
speed of $i$ upon entering the FEZ at time $t_{i}^{F}$ and let $u_{i}%
^{F}=u_{i}(t_{i}^{F})$ be the associated control (acceleration/deceleration).
Then, assuming for simplicity that a fixed control $u_{i}^{F}$ is maintained
throughout the FEZ, we have%
\[
F_{i}=\frac{\upsilon^{2}-(v_{i}^{F})^{2}}{2u_{i}^{F}},
\]
where $\upsilon\equiv v_{i}^{0}$ is the speed of $i$ when it reaches the CZ
after traveling a distance $F_{i}$. Clearly, the worst case in terms of the
maximal value of $F_{i}$, denoted by $\bar{F}_{i}$, arises when $k$
enters the CZ at minimal speed $v_{min}$ and $v_{i}^{F}=$ $v_{max}$, in which
case we must exert a minimal possible deceleration $u_B$, defined as a bound such that $u_{i,min} < u_B < 0$. Therefore,
\begin{equation}
\bar{F}_{i}=\frac{\upsilon^{2}-v_{max}^{2}}{2u_B} > 0. \label{Fworst}%
\end{equation}
On the other hand, recalling that $\tau=t_{i}^{0}$ is the time when $i$
reaches the CZ, the speed $\upsilon$ must also satisfy%
\begin{equation}
\tau-t_{k}^{0}=\frac{\upsilon-v_{max}}{u_B} > 0. \label{tauworst}%
\end{equation}
Thus, the length of the FEZ, $\bar{F}_{i}$, must be such that $(\tau
,\upsilon)\in\mathcal{F}_{i}$ subject to (\ref{Fworst})-(\ref{tauworst}). We
show next that under a sufficient condition on the system parameters $v_{min}%
$, $v_{max}$, $u_B$, and $\delta$, there exists a value of $\bar{F}_{i}$
which guarantees that $(\tau,\upsilon)\in\mathcal{F}_{i}.$

\begin{custompro}{1}
Suppose that
\begin{equation}
\frac{v_{min}-v_{max}}{u_B} \geqslant\frac{\delta}{v_{min}}
\label{parameter_condition}%
\end{equation}
holds. Then, the following FEZ length guarantees that $(\tau,\upsilon
)=(t_{k}^{0}+\frac{v_{min}-v_{max}}{u_B},v_{min})\in\mathcal{F}_{i}$:%
\begin{equation}
\bar{F}=\frac{v_{min}^{2}-v_{max}^{2}}{2u_B}. \label{F_upperbound}%
\end{equation}

\end{custompro}


\emph{Proof:} A necessary condition for the safety constraint (\ref{rearend}%
) to be satisfied throughout the CZ is that $s_{i}(\tau)\geqslant\delta$. This
is equivalent to $p_{k}(\tau)\geqslant\delta$, i.e., the distance traveled by
$k$ by the time CAV $i$ enters the CZ must be no less that the safety lower
bound $\delta$. The worst case arises when $v_{k}(t_{k}^{0})=v_{min}$ and
remains constant at least through $[t_{k}^{0},\tau]$. This implies that%
\begin{equation}
\tau-t_{k}^{0}\geqslant\frac{\delta}{v_{min}}.\label{timelagworst}%
\end{equation}
Moreover, observe that an upper bound for $\bar{F}_{i}$ in (\ref{Fworst}),
denoted by $\bar{F}$, occurs when $\upsilon=v_{min}$, so that
(\ref{F_upperbound}) holds. Then, (\ref{tauworst}) and (\ref{timelagworst})
imply (\ref{parameter_condition}). If this is satisfied, then $\tau=t_{k}%
^{0}+\frac{v_{min}-v_{max}}{u_B}$ is feasible, hence $(\tau,\upsilon
)=(t_{k}^{0}+\frac{v_{min}-v_{max}}{u_B},v_{min})\in\mathcal{F}_{i}$.
$\blacksquare$

Our analysis thus far has considered the case where the FEZ contains only CAV
$i$ and its preceding CAV $k$. This allows us to specify the upper bound
$\bar{F}$ in (\ref{F_upperbound}) for any such $i$. In general, however, there
may already be multiple CAVs in the FEZ at the time that a new CAV enters it.
We establish next that all such CAVs can be controlled to attain initial
conditions in their respective feasibility regions.

\begin{custompro}{2}
Let CAV $k$ enter the CZ when $N$ CAVs are in the
preceding FEZ, ordered so that $k<k_{0}<k_{1}<\cdots<k_{N}$ with associated
initial conditions when reaching the CZ $(\tau_{j},\upsilon_{j})$,
$j=0,\ldots,N$. Assume that (\ref{parameter_condition}) holds and the FEZ
length is given by (\ref{F_upperbound}). Then, $(\tau_{j},\upsilon_{j}%
)\in\mathcal{F}_{j}$ for all $j=0,\ldots,N$.
\end{custompro}
%

\emph{Proof:} From Proposition 1, setting $i=k_{0}$ we can attain $(\tau
_{0},\upsilon_{0})\in\mathcal{F}_{0}$. It follows that all information related
to $k_{0}$ is available to $k_{1}$ through the information set $Y_{k_{1}}(t)$
in (\ref{InfoSet}). Next, setting $k=k_{0}$ and $i=k_{1}$, we can again apply
Proposition 1 to attain $(\tau_{1},\upsilon_{1})\in\mathcal{F}_{1}$. This
iterative process is repeated over all $j=0,\ldots,N$. $\blacksquare$



\section{Simulation Examples}

The effectiveness of the proposed FEZ and associated control is illustrated
through simulation in MATLAB. For each direction, only one lane is considered.
The parameters used are: $L=400$ m, $S=30$ m, $\delta=10$ m, $v_{max}=15$ m/s,
$v_{min}=7$ m/s, $u_{i,max}=3$ m/s$^{2}$, $u_{i,min} = -5$ m/s$^2$ and $u_B=-2$ m/s$^{2}$, which
satisfy condition (\ref{parameter_condition}). Based on (\ref{F_upperbound}),
the length of the FEZ is set at $\bar{F}=44$ m. CAVs arrive at the FEZ based
on a random arrival process and any speed within $[v_{min},v_{max}]$. Here, we
assume a Poisson arrival process with rate $\lambda=1$ and the speeds are
uniformly distributed over $[7,15]$.

\begin{figure}[ptb]
\centering
\includegraphics[width = 3.2 in]{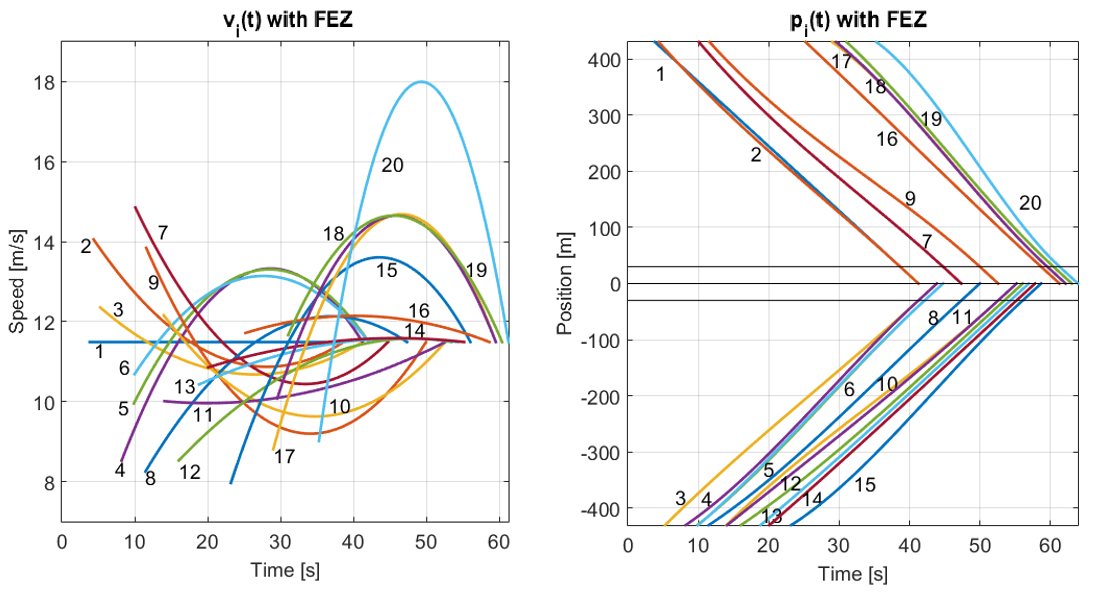} \caption{Speed $v_{i}(t)$ and
position $p_{i}(t)$ trajectories of the first 20 CAVs.}%
\label{pv}%
\end{figure}

We consider two cases: $(i)$ The FEZ is included preceding the CZ, and $(ii)$
No FEZ is included. The speed and position trajectories of the first 20
CAVs for the first case are shown in Fig. \ref{pv}. In the position
profiles, CAVs are separated into two groups: CAV positions shown above zero
are driving from east to west or from west to east, and those below zero are
driving from north to south or from south to north. These figures include
different instances from each of Cases 1), 2), or 3) in Section III.A
regarding the value of $t_{i}^{f}$. For example, CAV \#2 is assigned
$t_{2}^{f}=t_{1}^{f}$ and $v_{2}^{f}=v_{1}^{f}$, which corresponds to Case 1),
whereas CAV \#3 is assigned $t_{3}^{f}=t_{2}^{f}+\frac{S}{v_{3}^{f}}$ and
$v_{3}^{f}=v_{2}^{f}$, which corresponds to Case 3) with $v_{i}(t_{i-1}%
^{f})=v_{i-1}(t_{i-1}^{f})$.

\begin{figure}[ptb]
\centering
\includegraphics[width = 3.2 in]{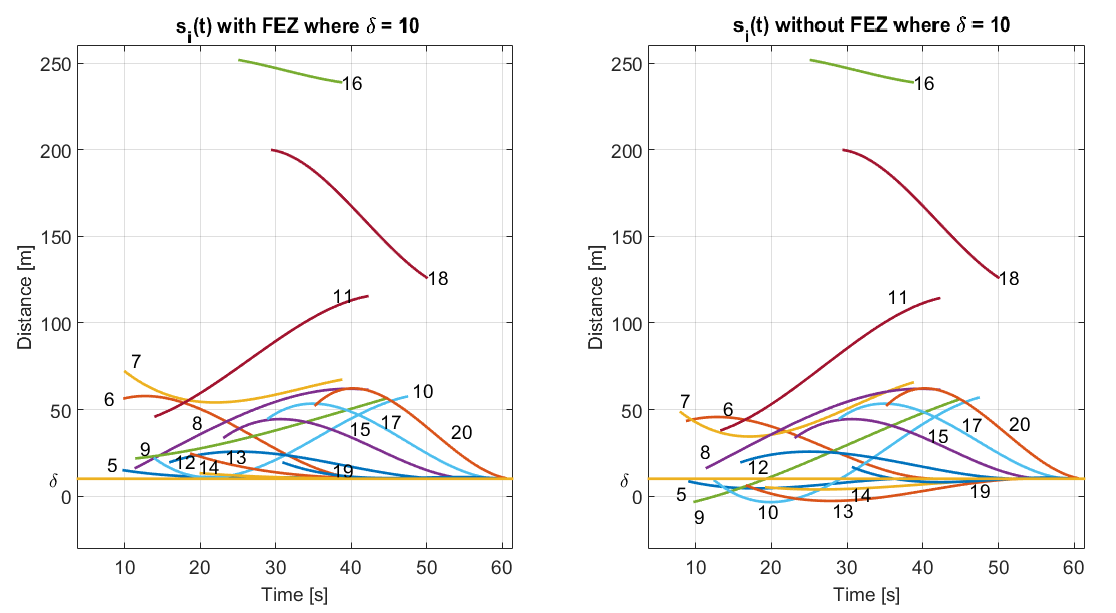} \caption{Distance
$s_{i}(t)$ trajectories of the first 20 CAVs.}%
\label{sim_sit}%
\end{figure}

To demonstrate the effectiveness of our feasibility enforcement control, we
examine the distance $s_{i}(t)$ between two consecutive CAVs for the first 20
CAVs as shown in Fig. \ref{sim_sit}. In this example, each of CAVs \#1, \#2,
\#3 and \#4 happens to be the first entering each of the four lanes
respectively, hence all $s_{i}(t)$, $i=1,\ldots,4$, in (\ref{rearend}) are
undefined. Regarding the different $s_{i}^{\ast}(\tau,\upsilon)$ cases arising
in Section IV, we observe that for CAV \#11, $t_{11}^{\ast}=t_{11}^{0}$, which
corresponds to Case 1; for CAV \#16, $t_{16}^{\ast}=t_{2}^{m}$, which
corresponds to Case 2; and for CAV \#7, $t_{7}^{\ast}=t_{1}\in(t_{7}^{0}%
,t_{1}^{m})$, which corresponds to Case 3. Without the FEZ (right side of Fig.
\ref{sim_sit}), we can see that CAVs \#5, \#9, \#10, \#13, \#14 and \#19 all
clearly violate the safety constraint (\ref{rearend}), i.e., $s_{i}%
(t)<\delta=10$ for at least some $t\in\lbrack t_{i}^{0},t_{i}^{m}]$. With the
FEZ included, these CAVs are capable of adjusting their speed and CZ entry
time to some feasible initial conditions and they all satisfy the safety
constraint, as clearly seen on the left side of Fig. \ref{sim_sit}. On the
other hand, given that CAV \#16 is on the same lane as CAV \#2 and that
CAV \#2 is the first one in that lane, there is no need for \#16 to make any
adjustment since it already has feasible initial conditions with respect to
the optimal control problem solved within the CZ.

\section{Conclusion}

Earlier work \cite{ZhangMalikopoulosCassandras2016} has established a
decentralized optimal control framework applied for coordinating online a
continuous flow of CAVs crossing two adjacent intersections in an urban area.
However, the feasibility of the optimal control solution \ depends on the
initial conditions of each CAV as it enters the \textquotedblleft control
zone\textquotedblright\ (CZ) of each intersection. We have shown that there
exists a feasibility region for each CAV in the space defined by its arrival
time and speed and this can be fully characterized in terms of information
known to CAV $i$ before it enters the CZ, which can be enforced through a
properly designed \textquotedblleft feasibility enforcement
zone\textquotedblright\ (FEZ) that precedes the CZ. Ensuring that optimal
control solutions are feasible paves the way for exploring more efficient
event-driven solutions, allow for different classes of CAVs with distinct
physical characteristics, and for alternative problem formulations that
exploit a potential trade-off between fuel consumption and congestion.

\bibliographystyle{IEEETran}
\bibliography{ACC2017_arXiv}

\end{document}